\documentclass[11pt,reqno]{amsart}
\usepackage[body={7in,9in},centering]{geometry}
\usepackage{amsmath,amssymb,amsthm,mathrsfs}
\usepackage{color,xcolor}
\definecolor{cobalt}{RGB}{61,99,181}
\usepackage[colorlinks,citecolor=cobalt,linkcolor=cobalt]{hyperref}
\usepackage{float}
\usepackage{exscale}
\usepackage{relsize}
\usepackage{graphicx}
\usepackage{tikz}
\usepackage{amsmath}

\numberwithin{equation}{section}

\date{\today}

\makeatletter

\newcommand{\Rmnum}[1]{\expandafter\@slowromancap\romannumeral #1@}
\makeatother

\usepackage{CJKutf8}

\title[A sufficient condition for the height function to be constant in $ I_g\times_\rho \mathbb{P}^n  $]{A sufficient condition for the height function to be constant in $ I_g\times_\rho \mathbb{P}^n  $}
\author[Kaijian Cao]{Kaijian Cao}
\address{
	\textsuperscript{1}
	College of Mathematics and Statistics, Chongqing University, Chongqing, 401331, P. R. China}
\email{20175785@cqu.edu.cn}

\keywords{warped product space; Omori-Yau maximum principle; Newton tensor; 2-mean curvature.}

\begin{document}
	\begin{CJK}{UTF8}{gbsn}
	\maketitle

	\section{abstract}
	This paper makes some modifications to the warped product space. Based on Alias,Impera and Rigoli, a warping function is added to the warped product space. This new function affects the Riemannian metric of the warped product space. In this new warped product space, we continue to discuss the sufficient condition for calculating the height of the immersed surface.
	
	\section{Preparation}
	\subsection{ warped product space $ I_g\times _\rho \mathbb{P}^n  $}
	\leavevmode \\ \\
	when we talk about $ I_g\times _\rho \mathbb{P}^n  $，among $ I\subset \mathbb{R} $ is an open interval，$ \mathbb{P}^n $is a complete n-dimensional Riemannian manifold,
	$$ \rho:I\rightarrow \mathbb{R}^+ $$
	$$ g:\mathbb{P}^n \rightarrow \mathbb{R}^+$$
	Both are positive warping functions.
	
	The Riemannian metric for such a warped product space  $ I_g\times _\rho \mathbb{P}^n $is：
	$$<,>=g^2(\pi_{\mathbb{P}})\pi_I^\star(dt^2)+\rho^2(\pi_I)\pi_{\mathbb{P}}^\star(<,>_{\mathbb{P}}),$$
	Here $ \pi_I:I_g\times _\rho \mathbb{P}^n \rightarrow I $，$\pi_{\mathbb{P}}:I_g\times _\rho \mathbb{P}^n \rightarrow \mathbb{P} $ represents the natural projection to its two subspaces of $  I_g\times _\rho \mathbb{P}^n$respectively.
	
	\subsection{Some concepts}
	\leavevmode \\ \\
	\textbf{Definition2.2.1}:Let$ f:\Sigma^n\rightarrow \bar{M}^{n+1} $be an isometric immersion from $ \Sigma^n $ to $ \bar{M}^{n+1} $,and$ \forall p\in \Sigma,X,Y\in T_p \Sigma$ we have：
	$$ A: T_p \Sigma\rightarrow T_p \Sigma$$
	$$ <AX,Y>=-<\bar{\nabla}_XN,Y>    ,$$
	then $ A $ is the shape operator with respect to unit normal field  $ N $， $\bar{\nabla}  $ is Riemannian connection of $ \bar{M}^{n+1} $.
	
	Let $ \lambda_1,\dots,\lambda_n $ be the eigenvalue of $ A $，记$ e_1,\dots,e_n $are the corresponding orthonormal eigenvectors. The symmetric function is given by the eigenvalues as follows：
	\begin{align}
	\nonumber	S_0&=1 \\
	\nonumber	S_k&=\Sigma_{i_1<\dots i_k}\lambda_{i_1}\lambda_{i_2}\dots\lambda_{i_k},\quad 1\le k\le n\\
	\nonumber	S_k&=0, \quad  k>n.	
	\end{align}

\noindent \textbf{Definition2.2.2}:Under the above assumptions，$ k- $mean curvature $ H_k $ is defined as follows定义如下：
$$ H_k=\dfrac{S_k}{C_n^k}  .  $$

\noindent\textbf{Definition2.2.3}:Under the above assumptions， Nweton tensor is defined as follows:
\begin{align}
\nonumber	P_0&=I,\\
\nonumber	P_k&=S_kI-AP_{k-1},\quad k=1,2,\dots,n.
\end{align}

It's easy to see that $ P_k=S_kI-S_{k-1}A+\dots (-1)^kA^k。 $ and obviously $   P_k$ is a symmetric linear map.

\noindent\textbf{Definiton2.2.4}:
\begin{align}
\nonumber	L_k&:C^\infty(\Sigma)\rightarrow C^\infty(\Sigma)\\
\nonumber	L_k&=Tr(P_k\circ hess).
\end{align}

\subsection{Some properties}
	\leavevmode \\ \\
\noindent\textbf{Proposition2.3.1}:$ \forall 1\le k\le n $, if $ H_1,H_2,\dots,H_n>0 $, then:
$$   H_{k-1}H_{k+1}\le H_k^2.     $$

\noindent\textbf{Proposition2.3.2}:$ \forall 1\le k\le n $, if $ H_1,H_2,\dots,H_n>0 $, then:
\begin{align}
\nonumber	(1)H_1^2&\ge H_2;\\
\nonumber	(2)H_1&\ge H_2^{\frac{1}{2}}\ge H_3^{\frac{1}{3}}\ge \dots \ge H_n^{\frac{1}{n}}.
\end{align}

\noindent\textbf{Proposition2.3.3}: If $ H_2>0 $, then $ L_1=Tr(P_1\circ hess) $ is elliptic operator.

\section{A generation of Omori-Yau maximum principle}

\subsection{Omori-Yau Maximum Principle}
	\leavevmode \\ \\
\noindent\textbf{Definition3.1.1}:Let $ \Sigma $ be n-dimensional Riemannian manifold, if $\forall u\in C^2(\Sigma),u^\star=\sup_{\Sigma}u<+\infty,$  there exists $\{p_j\}\subset\Sigma$ such that：
\begin{align}
\nonumber	(1)&u(p_j)>u^\star-\frac{1}{j}\\
\nonumber	(2)&||\nabla u(p_j)||<\frac{1}{j}\\
\nonumber	(3)&\Delta u(p_j)<\frac{1}{j}
\end{align}
then we call Omori-Yau maximum principle hold for $ \Delta $ in $ \Sigma $.

Equivalently, we can also describe another equivalent description，if $\forall u\in C^2(\Sigma),u^\star=\inf_{\Sigma}u>-\infty,$ there exists $\{p_j\}\subset\Sigma$ such that：
\begin{align}
	\nonumber	(1)&u(p_j)<u_\star+\frac{1}{j}\\
	\nonumber	(2)&||\nabla u(p_j)||<\frac{1}{j}\\
	\nonumber	(3)&\Delta u(p_j)>-\frac{1}{j}
\end{align}
then Omori-Yau maximum principle hold for $ \Delta $ in $ \Sigma $.
\leavevmode \\

\noindent\textbf{Definition3.1.2}:Let $ \Sigma $ ba a n-dimensional Riemannian manifold，$ L=Tr(P\circ hess) $ is a semidefinite elliptic operator，$ P:T\Sigma \rightarrow  T\Sigma$ is a semidefinite symmetric operator satisfying $ \sup_{\Sigma}Tr(P)<+\infty $. If $\forall u\in C^2(\Sigma),u^\star=\sup_{\Sigma}u<+\infty,$ there exists $\{p_j\}\subset\Sigma$such that：
\begin{align}
	\nonumber	(1)&u(p_j)>u^\star-\frac{1}{j}\\
	\nonumber	(2)&||\nabla u(p_j)||<\frac{1}{j}\\
	\nonumber	(3)&Lu(p_j)<\frac{1}{j}
\end{align}
then下Omori-Yau maximum principle holds for  $ L $ in $ \Sigma $.

\subsection{A sufficient condition for 下Omori-Yau maximum principle under $ L $}
	\leavevmode \\ \\
\noindent\textbf{Theorem3.2.1}:Let $ (\Sigma,<>) $ be a Riemannian manifold, $ L=Tr(P\circ hess) $ is a semidefinite elliptic operator，$ P:T\Sigma \rightarrow T\Sigma $ is a semidefinite symmetric operator satisfying $ \sup_{\Sigma}Tr(P)<+\infty $. If there exists non-negative function $ \gamma \in C^2(\Sigma) $ such that:
\begin{align}
\nonumber	(1)&\gamma(p)\rightarrow +\infty,\quad p\rightarrow \infty \\
\nonumber	(2)&\exists A>0,s.t.||\nabla\gamma||\le A\gamma^{\frac{1}{2}},\quad \text{out side a compact set}\\
\nonumber	(3)&\exists B>0,s.t.L\gamma\le B\gamma^{\frac{1}{2}}G^{\frac{1}{2}}(\gamma\frac{1}{2}),\quad \text{out side a compact set}.
\end{align}
$ G $ is a smooth function in $ [0,+\infty) $ satisfying:
\begin{align}
\nonumber	(i)&G(0)>0\\
\nonumber	(ii)&G'(t)\ge 0,\quad \text{in }[0,+\infty)\\
\nonumber	(iii)&G^{-\frac{1}{2}}(t)\notin L^{1}[0,+\infty)\\
\nonumber	(iv)&\limsup_{t\rightarrow +\infty}\dfrac{tG(t\frac{1}{2})}{G(t)}<+\infty.
\end{align}

then Omori-Yau maximum principle holds for  $ L $ in $ \Sigma $.
\leavevmode \\ \\

Using theorem 3.2.1， we want  Omori-Yau maximum principle holds in $ \Sigma $, so naturally wo need to find appropriate $ \gamma $ and $G  $.

We give some notations for convenience, Let$ (\Sigma,<>) $ be a complete non-compact Riemannian manifold, $ o\in \Sigma $ is a fixed point, $ r(p) $ is the distance function starting from $ o $ to $p$, $\gamma(p)=r^{2}(p)$.
\leavevmode \\
\noindent\textbf{Theorem3.2.2}:Let $ (\Sigma,<>) $ be a complete non-compact Riemannian manifold, satisying $ K_{\Sigma}^{rad}\ge -G(r) $, here $ K_{\Sigma}^{rad}  $ represents radial sectional curvature(sectional curvature including $ \nabla r $), $ G $ is the function which satisfies previous theorems and smooth at 0.

\subsection{Omori-Yau maximum principle in $ I_g\times _\rho \mathbb{P}^n  $}
	\leavevmode \\ \\
Here we consider $ \mathbb{P}^n $ to be a complete non-compact Riemannian manifold，$ o\in  \mathbb{P}^n  $ is the fixed point，$ \hat{r} $ represents distance fuction in $ \mathbb{P}^n $ starting from fixed point $ o $.

Following the previous setting for radial sectional curvature, we set:
$$  K_{\mathbb{P}}^{rad}\ge -G(\hat{r}) ,   $$ here $ G $ is the function which satisfies previous theorems and smooth at 0.

$ f $ is a properly immersed hypersurface, $ f(\Sigma)\subset [t_1,t_2]\times \mathbb{P}^n .$

To obtain the Omori-Yau maximum principle under some operator, we just need to use the previous theorem 3.2.1 to verify that our selected $\gamma $satisfies the following three conditions:
\begin{align}
	\nonumber	(1)&\gamma(p)\rightarrow +\infty,\quad p\rightarrow \infty \\
	\nonumber	(2)&\exists A>0,s.t.||\nabla\gamma||\le A\gamma^{\frac{1}{2}},\quad \text{outside a compact set}\\
	\nonumber	(3)&\exists B>0,s.t.L\gamma\le B\gamma^{\frac{1}{2}}G^{\frac{1}{2}}(\gamma\frac{1}{2}),\quad \text{outside a compact set}.
\end{align}
其中$ G $ is a smooth function in $ [0,+\infty) $ satisfying:
\begin{align}
	\nonumber	(i)&G(0)>0\\
	\nonumber	(ii)&G'(t)\ge 0,\quad \text{in }[0,+\infty)\\
	\nonumber	(iii)&G^{-\frac{1}{2}}(t)\notin L^{1}[0,+\infty)\\
	\nonumber	(iv)&\limsup_{t\rightarrow +\infty}\dfrac{tG(t\frac{1}{2})}{G(t)}<+\infty.
\end{align}

(1) Because $ f $ is proper, so $ f(p)\rightarrow \infty(p\rightarrow \infty) $, and $f(\Sigma)\subset[t_1,t_2]\times \mathbb{P}$, thus $ \hat{r}^2(\pi_{\mathbb{P}} \circ f(p))\rightarrow+\infty(p\rightarrow \infty) $, that is $ \gamma(p)=\tilde{\gamma}(f(p))= \hat{\gamma}(\pi_{\mathbb{P}}\circ f(p))=\hat{r}^2(\pi_{\mathbb{P}} \circ f(p))\rightarrow +\infty(p\rightarrow \infty).$\\

(2) We use $ \tilde{\nabla},\hat{\nabla},\nabla $ represent Riemannian connections in$ I_g\times _\rho \mathbb{P}^n ,\mathbb{P},\Sigma $ respectively, $ \gamma(p)=\tilde{\gamma}\circ f(p), $:
$$ \tilde{\nabla}\tilde{\gamma }=\nabla \gamma+<\tilde{\nabla}\tilde{\gamma },N>N $$ $ N $ is unit normal field in immersed hypersurface $ f $.

And $ \tilde{\gamma}\circ f(p)=\hat{\gamma}(\pi_{\mathbb{P}}\circ f(p)) $ has nothing to do with $ I $, so:
$$    <\tilde{\nabla}\tilde{\gamma },T>=0    $$ here $ T $ is pullback  of $ \dfrac{\partial}{\partial t} $ from $ TI $ to $ T(I_g\times _\rho \mathbb{P}^n) $.

From the gradient vector field transformation relation:
$$  <\tilde{\nabla}\tilde{\gamma },V>=<\hat{\nabla}\hat{\gamma},V>_{\mathbb{P}}      $$ here $ V $ is pullback from$ TP $ to $ T(I_g\times _\rho \mathbb{P}^n) $.

And according to the corresponding relationship of measurement:
$$  <\tilde{\nabla}\tilde{\gamma },V>=\rho^2<\tilde{\nabla}\tilde{\gamma },V>_{\mathbb{P}}        ,$$ so $ \tilde{\nabla}\tilde{\gamma }=\dfrac{1}{\rho^2}\hat{\nabla}\hat{\gamma}=\dfrac{2\hat{r}}{\rho^2}\hat{\nabla}\hat{r}, $ thus $ ||\hat{\nabla}\hat{r}||=\rho||\hat{\nabla}\hat{r}||_{\mathbb{P}} =\rho\ge \min_{[t_1,t_2]}\rho>0,$ that is $ ||\nabla\gamma||\le||\tilde{\nabla}\tilde{\gamma }||=\dfrac{2\gamma^{\frac{1}{2}}}{\rho}\le c\gamma^{\frac{1}{2}} .$\\

(3)
$$ Hess\gamma(X,X)=Hess\tilde{\gamma}(X,X)+<\tilde{\nabla}\tilde{\gamma },N><AX,X>,\quad X\in T\Sigma $$

$$ \tilde{\nabla}_T\tilde{\nabla}\tilde{\gamma}=-\mathcal{H}\tilde{\nabla}\tilde{\gamma}+\dfrac{1}{\rho^4g^2}<\hat{\nabla}\hat{\gamma},g\hat{\nabla}g>T (\mathcal{H}=\frac{\rho'}{\rho}) $$

$ Hess\tilde{\gamma}(T,T)=\dfrac{g}{\rho^2}<\tilde{\nabla}\tilde{\gamma},\hat{\nabla}g>=\dfrac{g}{\rho^4}<\hat{\nabla}\hat{\gamma},\hat{\nabla}g>. $

$ X\in T\Sigma $ is decomposed as follows:
$$     X=X^\star+\dfrac{1}{g^2}<X,T>T  $$ here $ X=d\pi_{\mathbb{P}}(X). $

then:
$$  Hess\tilde{\gamma}(X,X)=Hess\tilde{\gamma}(X^\star,X^\star)+\dfrac{2<X,T>}{g^2}Hess\tilde{\gamma}(X^\star,T)+\dfrac{Hess\tilde{\gamma}(T,T)}{g^4}<X,T>^2      $$
and:
$$  Hess\tilde{\gamma}(X^\star,T)=-\mathcal{H}<\nabla\gamma,X> .   $$

we use:
$$ \tilde{\nabla}_{X^\star}\tilde{\nabla}\tilde{\gamma}=\dfrac{1}{\rho^2}\hat{\nabla}_{X^\star}\hat{\nabla}\hat{\gamma}-\dfrac{\rho'}{\rho^3g^2}<\hat{\nabla}\hat{\gamma},X^\star>T$$

and:
 \begin{align}
\nonumber 	Hess\tilde{\gamma}(X^\star,X^\star)&=\dfrac{1}{\rho^2}<\hat{\nabla}_{X^\star}\hat{\nabla}\hat{\gamma},X^\star>\\
\nonumber &=<\hat{\nabla}_{X^\star}\hat{\nabla}\hat{\gamma},X^\star>_{\mathbb{P}}\\
\nonumber &=Hess\hat{\gamma}(X^\star,X^\star). 
 \end{align}

So we simplify above by combining:
\begin{align}
\nonumber	Hess\gamma(X,X)&=Hess\tilde{\gamma}(X,X)+<\tilde{\nabla}\tilde{\gamma},N><AX,X>\\
\nonumber	&=Hess\tilde{\gamma}(X^\star,X^\star)+\dfrac{2<X,T>}{g^2}Hess\tilde{\gamma}(X^\star,T)+\dfrac{Hess\tilde{\gamma}(T,T)}{g^4}<X,T>^2+<,\tilde{\nabla}\tilde{\gamma},N><AX,X>\\
\nonumber	&=Hess\hat{\gamma}(X^\star,X^\star)-\dfrac{2<X,T>}{g^2}\mathcal{H}<\nabla,X>+\dfrac{<\hat{\nabla}\hat{\gamma},\hat{\nabla}g>}{g^3\rho^4}<X,T>^2+<\tilde{\nabla}\tilde{\gamma},N><AX,X>.
\end{align}

Computing $ Hess\hat{\gamma}(X^\star,X^\star) $ by Hessian comparison theorem:
$$  Hess\hat{\gamma}(X^\star,X^\star)\le(\gamma G(\gamma^{\frac{1}{2}}))^{\frac{1}{2}} ||X||^2 ,$$

and $ |\dfrac{2<X,T>}{g^2}<\nabla\gamma,X>|\mathcal{H}\le\dfrac{2||\nabla||||X||^2}{g}\mathcal{H}\dfrac{c}{g}\gamma^{\frac{1}{2}}||X||^2. $

And we give following restrictions on $ g $:
$$   g\ge 1,<,\hat{\nabla}\hat{\gamma},\hat{\nabla}g>=0     $$

then $Hess\gamma(X,X)\le c (\gamma G(\gamma^{\frac{1}{2}}))^{\frac{1}{2}}||X||^2+<\tilde{\nabla}\tilde{\gamma},N><AX,X>$.

The discussion is divided into the following two operators:

a. When the operator is $ \Delta $, suppose $ \sup_{\Sigma}|H|<+\infty $:
$$ \Delta \gamma=cn(\gamma G(\gamma^{\frac{1}{2}}))^{\frac{1}{2}}+nH<\tilde{\nabla}\tilde{\gamma},N>       ,$$
$ |H<\tilde{\nabla}\tilde{\gamma},N>|\le \sup_{\Sigma}|H|c\gamma^{\frac{1}{2}}\le C \gamma^{\frac{1}{2}}\le C(\gamma^{\frac{1}{2}}))^{\frac{1}{2}} \quad \gamma\text{is suffcient large}$

So there must be some compact set, and outside that compact set there is $ \delta \gamma \le C(\gamma^{\frac{1}{2}}))^{\frac{1}{2}} $\\

b.When the operator is $L$, suppose $\sup_{\Sigma}||A||^2< +\infty$:
$$ |<\tilde{\nabla}\tilde{\gamma},N><AX,X>|\le ||\tilde{\nabla}\tilde{\gamma}||||A||||X||^2\le C(\gamma^{\frac{1}{2}}))^{\frac{1}{2}}||X||^2  \quad \gamma\text{sufficient large}$$

\leavevmode \\ \\

The following theorem is obtained by summarizing the above discussion:

\noindent\textbf{Theorem3.3.1}:Let $ \mathbb{P}^n $ be a complete non-compact Riemannian manifold, its radial sectional curvature  satisfies $ K_{\mathbb{P}}^{rad}\ge-G(\hat{r}) $，$ f:\Sigma\rightarrow I_g\times_\rho \mathbb{P} $ is  properly immersed hypersurface, and$ f(\Sigma)\subset [t_1,t_2]\times\mathbb{P}^n,g\ge 1,<\hat{\nabla}\hat{\gamma},\hat{\nabla}g>=0 $， if $ \sup_{\Sigma}|H|<+\infty $, then  Omori-Yau maximum principle holds for $ \Delta $ in $ \Sigma $.。

\noindent\textbf{Theorem3.3.2}:Let $ \mathbb{P}^n $ be a complete non-compact Riemannian manifold, its radial sectional curvature  satisfies $ K_{\mathbb{P}}^{rad}\ge-G(\hat{r}) $，$ f:\Sigma\rightarrow I_g\times_\rho \mathbb{P} $ is a  properly immersed hypersurface, and$ f(\Sigma)\subset [t_1,t_2]\times\mathbb{P}^n,g\ge 1,<\hat{\nabla}\hat{\gamma},\hat{\nabla}g>=0 $， if $ \sup_{\Sigma}||A||<+\infty $，then for all semidefinite elliptic operator $ L=Tr(P\circ hess) ,$ $\sup_{\Sigma}Tr(P)<+\infty$,  Omori-Yau maximum principle holds for $ L $ in $ \Sigma $.

\section{Conclusion}
\noindent\textbf{Proposition4.1.1}:Let $ f:\Sigma \rightarrow I_g\times_\rho \mathbb{P} $ be an isometric immersion, $ h=\pi_I\circ f $ is the height function, $ \sigma(t)=\int_{t_0}^{t}\rho(u)du,<P_k\nabla h,(\hat{\nabla}\frac{1}{g^2})^T>=0 $,  then：
\begin{align}
\nonumber	L_kh&=\dfrac{\rho'}{g^2\rho}(c_kH_k-g^2<P_k\nabla h,\nabla h>)+\dfrac{c_k\theta H_{k+1}}{g^2}\\
\nonumber	L_k\sigma(h)&=\dfrac{\rho'}{g^2}c_kH_k+\dfrac{c_k\rho \theta H_{k+1}}{g^2}
\end{align}
其中$c_k=(n-k)C_n^k=(k+1)C_n^{k+1},\theta=<N,T>.$

\begin{proof}
	$ \bar{\nabla}\pi_I=\dfrac{T}{g^2}, $ then $\nabla h=(\bar{\nabla}\pi_I)^T=\dfrac{1}{g^2}(T-\theta N)$, so $ T^T=T-\theta N=g^2\nabla h $.
	
	Using $ \bar{\nabla}_VT=\bar{\nabla}_TV=\dfrac{1}{\rho^2g^2}<V,g\hat{\nabla}g>T+(Tlin\rho)V,\quad V\in TP. $
	
	$\bar{\nabla}_TT=-\dfrac{g\hat{\nabla}g}{\rho^2}$
	
	So
	\begin{align}
\nonumber		\bar{\nabla}_XT&=\bar{\nabla}_{X^\star+\frac{1}{g^2}<X,T>T}T\\
\nonumber		&=-\dfrac{1}{\rho^2g}<X,T>\hat{\nabla}g+\dfrac{1}{\rho^2g^2}<X^\star,g\hat{\nabla}g>T+(Tln\rho)X^\star
	\end{align}
	here $ X\in T(I_g\times_\rho \mathbb{P}). $
	
	Now we discuss $ X\in T\Sigma $:
	\begin{align}
\nonumber		\bar{\nabla}_X\nabla h&=\bar{\nabla}_X(\dfrac{1}{g^2}(T-\theta N))\\
\nonumber		&=X(\dfrac{1}{g^2})(T-\theta N)+\dfrac{1}{g^2}\bar{\nabla}_XT-\dfrac{X(\theta)N}{g^2}+\dfrac{\theta}{g^2}AX.
	\end{align}

So we have:
\begin{align}
\nonumber	hess h(X)&=\nabla_X\nabla h=(\bar{\nabla}_X\nabla h)^T\\
\nonumber	&=g^2X(\dfrac{1}{g^2})\nabla h-\dfrac{1}{\rho^2g}<X,\nabla h>(\hat{\nabla} g)^T+\dfrac{1}{\rho^2g}<X,(\hat{\nabla} g)^T>\nabla h+\dfrac{\rho'}{\rho g^2}X^\star+\dfrac{\theta}{g^2}AX.
\end{align}	

	So we can calculate $ L_kh ,L_k\sigma(h)$ as follows:
	\begin{align}
\nonumber		L_kh&=Tr(P_k\circ hess h)=\sum_{i}<P_k\circ hess h(e_i),e_i>\\
\nonumber		&=g^2<P_k\nabla h,(\hat{\nabla}\frac{1}{g^2})^T>+\dfrac{\theta}{g^2}Tr(P_kA)+\dfrac{\rho'}{\rho g^2}Tr(P_k)-\dfrac{\rho'}{\rho}<P_k\nabla h,\nabla h>\\
\nonumber		&=\dfrac{\rho'}{g^2\rho}(c_kH_k-g^2<P_k\nabla h,\nabla h>)+\dfrac{c_k\theta H_{k+1}}{g^2}.
	\end{align}

	\begin{align}
\nonumber		L_k\sigma(h)&=Tr(P_k\circ hess \sigma(h))\\
\nonumber		&=\sum_{i}<e_i(\rho)\nabla h,P_ke_i>+\rho \sum_{i}<\nabla_{e_i}\nabla h,P_ke_i>\\
\nonumber		&=<P_k\nabla h,\rho'\nabla h>+\rho L_kh\\
\nonumber		&=\dfrac{\rho'}{g^2}c_kH_k+\dfrac{c_k\rho\theta H_{k+1}}{g^2}.
	\end{align}
\end{proof}

\noindent\textbf{Proposition4.1.2}:Let $f:\Sigma\rightarrow I_g\times_\rho \mathbb{P}  $ be a hypersurface whose mean curvature is not zero everywhere，$ f(\Sigma)\subset [t_1,t_2]\times\mathbb{P} $, if $ \mathcal{H}\ge 0,<\nabla h,(\hat{\nabla}\frac{1}{g^2})^T>=0 $, and the angle function $ \theta $ doesn't change sign， then after choosing a suitable direction in $ \Sigma $ such that $ H_1>0 $, we know that Omori-Yau maximum principle holds for $ \Delta $ in $\Sigma$, then:

(1) if $ \theta \le 0, $ then $\rho'\ge 0$\\
(2) if $ \theta \ge 0, $ then $\rho'\le 0$

\begin{proof}
	Because Omori-Yau maximum principle holds for $ \Delta $ in  $\Sigma$, We can find a list of points $ \{p_j\}\subset \Sigma $ satisfying:
	\begin{align}
\nonumber		\lim_{j\rightarrow+\infty}h(p_j)&=h_\star=\inf_{\Sigma}h\\
\nonumber		||\nabla h(p_j)||^2&=\dfrac{1}{g^4}(g^2-\theta^2)<(\dfrac{1}{j})^2\\
\nonumber		\Delta h(p_j)&=\dfrac{\rho'}{\rho g^2}(n-g^2||\nabla h(p_j)||^2)+\dfrac{n\theta H_1}{g^2}>-\dfrac{1}{j}
	\end{align}
	
	then $-\dfrac{n\theta H_1}{g^2}<\dfrac{1}{j}+\dfrac{\rho'}{\rho g^2}(n-g^2||\nabla h(p_j)||^2)$
	
	Analogously, for $ h^\star $ we still have $ \{q_j\}\subset \Sigma $ satisfying:
		\begin{align}
\nonumber		\lim_{j\rightarrow+\infty}h(q_j)&=h^\star=\sup_{\Sigma}h\\
\nonumber		||\nabla h(q_j)||^2&=\dfrac{1}{g^4}(g^2-\theta^2)<(\dfrac{1}{j})^2\\
\nonumber		\Delta h(q_j)&=\dfrac{\rho'}{\rho g^2}(n-g^2||\nabla h(q_j)||^2)+\dfrac{n\theta H_1}{g^2}<\dfrac{1}{j}
	\end{align}

 then $-\dfrac{n\theta H_1}{g^2}>-\dfrac{1}{j}+\dfrac{\rho'}{\rho g^2}(n-g^2||\nabla h(q_j)||^2)$

	There are two conditions to consider:
	
	(1)$ \theta\le 0 $:
	
	Then $-\theta(p_j)\ge 0,H_1(p_j)>0,$ using $-\dfrac{n\theta H_1}{g^2}<\dfrac{1}{j}+\dfrac{\rho'}{\rho g^2}(n-g^2||\nabla h(p_j)||^2)$， let $ j\rightarrow \infty $ we have $\rho'\ge 0.$
	
	(2)$ \theta\ge 0 $:	
	
		Then $\theta(q_j)\ge 0,H_1(q_j)>0,$ using $-\dfrac{n\theta H_1}{g^2}>-\dfrac{1}{j}+\dfrac{\rho'}{\rho g^2}(n-g^2||\nabla h(q_j)||^2)$, let $ j\rightarrow \infty $ we have $\rho'\le 0.$
		
\end{proof}

\subsection{A sufficient condition for the height function to be constant in $ I_g\times_\rho \mathbb{P}^n  $}
	\leavevmode \\ \\

Define $ \xi_1=(n-1)\dfrac{\rho'}{\rho}\Delta-\theta L_1 $, so $\xi_1\sigma(h)=\dfrac{n(n-1)\rho}{g^2}(\dfrac{\rho'^2}{\rho^2}-\theta^2H_2)  $, then $\xi_1=Tr(\zeta_1\circ hess)  $, here$ \zeta_1=(n-1)\dfrac{\rho'}{\rho}I-\theta P_1. $

\noindent\textbf{Theorem4.2.1}:Let $f:\Sigma\rightarrow I_g\times_\rho \mathbb{P}  $ be an compact immersed hypersurface  and $ H_2>0 $, if $ \mathcal{H}\ge 0 $, $ \theta $ doesn't change sign, $ <P_1\nabla h,(\hat{\nabla}\dfrac{1}{g^2})^T>=<\nabla h,()\hat{\nabla}\dfrac{1}{g^2}^T>=0 $, and for maximum and minimum points $ p_{max},p_{min} $ of $ h $ we have $ g(p_{max})=g(p_{min}),\theta(p_{max})=\theta_{max},\theta(p_{min})=\theta_{min} $. Then $ h $ is a constant function.

\begin{proof}
	(1)$ \theta\le 0,\rho'\ge 0 $:
	
	$ \zeta_1e_i=[(n-1)\dfrac{\rho'}{\rho}-\theta(S_1-\lambda_i)]e_i $， obviously this is a semi-definite operator，$ ||\nabla h(p_{max})||=||\nabla h(p_{min})||=0 $, so:
	\begin{align}
\nonumber		\theta(p_{max})&=-g(p_{max})\\
\nonumber		\theta(p_{min})&=-g(p_{min})
	\end{align}
$ 	\theta(p_{max})=	\theta(p_{min}) $.

And $ \sigma $ is monotonically increasing:
\begin{align}
\nonumber	(\sigma\circ h)^\star&=\max_{\Sigma}\sigma\circ h=\sigma(h^\star)=\sigma(h(p_{max}))\\
\nonumber	(\sigma\circ h)_\star&=\min_{\Sigma}\sigma\circ h=\sigma(h_\star)=\sigma(h(p_{min}))
\end{align}
	
	Because of the extreme point property:
	\begin{align}
\nonumber		Hess\sigma(h)(p_{max})&\le 0\\
\nonumber				Hess\sigma(h)(p_{min})&\ge 0
	\end{align}
	
So:
\begin{align}
\nonumber	\xi_1\sigma(h(p_{max}))&=\dfrac{n(n-1)\rho}{g^2}(h^\star)(\dfrac{\rho'^2}{\rho^2}(h^\star)-\theta^2(p_{max})H_2)\le 0\\
\nonumber		\xi_1\sigma(h(p_{min}))&=\dfrac{n(n-1)\rho}{g^2}(h_\star)(\dfrac{\rho'^2}{\rho^2}(h_\star)-\theta^2(p_{min})H_2)\ge 0
\end{align}
 
 Thus $ \mathcal{H}(h_\star)=\mathcal{H}(h^\star) $，$ \mathcal{H} $ is a  constant function.
 
 So we have：
 \begin{align}
\nonumber 	L_1(\sigma\circ h)&=\dfrac{n(n-1)\rho}{g^2}(\dfrac{\rho'}{\rho}H_1+\theta H_2)\\
\nonumber 	 &\ge \dfrac{n(n-1)\rho}{g^2}(\dfrac{\rho'}{\rho}H_1+\theta_{min} H_2)\\
\nonumber 	 &=\dfrac{n(n-1)\rho}{g^2}H_2^\frac{1}{2}g(p_{min})(H_1-H_2^\frac{1}{2})\ge 0
 \end{align}

From the extreme value principle of elliptic operators, $ h $ is a constant function.

(2)$ \theta\ge 0,\rho'\le 0 $ has the same discussion, we only need to consider $ -\zeta_1 $.

\end{proof}

\noindent\textbf{Theorem4.2.2}: Let $f:\Sigma\rightarrow I_g\times_\rho \mathbb{P}  $ be an compact immersed hypersurface  and $ H_2>0 $，$ K_{\Sigma}^{rad}\ge -G(r) $,$ K_{\Sigma}^{rad} $ represents radial sectional curvature, $ G $ is a smooth function that satisfies  previous theorems. If $ \sup_{\Sigma}|H_1|<+\infty,f(\Sigma)\subset [t_1,t_2]\times\mathbb{P},\mathcal{H}'>0, \theta$ doesn't change sign and attains its maximum and minimum, and there exists $ p_{max},p_{min}\in \Sigma $ such that $ \theta(p_{max})=\theta_{max},\theta_{p_{min}}=\theta_{min},\lim_{p\rightarrow \infty}g(p)=1,<P_1\nabla h,(\hat{\nabla}\dfrac{1}{g^2})^T>=<\nabla h,(\hat{\nabla}\dfrac{1}{g^2})^T>=0 $, then $ h $ is a constant function.

\begin{proof}
	(1)$\theta\le 0,\rho'\ge 0 $: 
	
	$ \zeta_1 $ is a semidefinite operator $ Tr(\zeta_1)=n(n-1)\dfrac{\rho'}{\rho}-n(n-1)\theta H_1 \le n(n-1)(\dfrac{\rho'}{\rho}(h_\star)-\theta_{min}H_1^\star)<+\infty$,$ H_1^\star=\sup_{\Sigma}H_1<+\infty .$
	
	From theorem3.2.2, Omori-Yau maximum principle holds for $ \xi_1 $ in $\Sigma$.
	
	So we can find a list of points $\{p_j\}\subset \Sigma$ such that:
	 	\begin{align}
\nonumber	 	\lim_{j\rightarrow+\infty}(\sigma\circ h)(p_j)&=\sup_{\Sigma}(\sigma\circ h)=\sigma(h^\star)<+\infty\\
\nonumber	 	||\nabla (\sigma\circ h)(p_j)||&=\rho(h(p_j))||\nabla h(p_j)||<\dfrac{1}{j}\\
\nonumber	 	\xi_1 (\sigma\circ h)(p_j)&<\dfrac{1}{j}
	 \end{align}
 Let $ j\rightarrow +\infty $ we have $ \mathcal{H}(h^\star)\le H_2 .$
 
Equivalently from Omori-Yau maximum principle we have a list of pioints $\{q_j\}\subset \Sigma$ such that:
 \begin{align}
\nonumber 	\lim_{j\rightarrow+\infty}(\sigma\circ h)(q_j)&=\inf_{\Sigma}(\sigma\circ h)=\sigma(h_\star)\\
\nonumber 	||\nabla (\sigma\circ h)(q_j)||&=\rho(h(q_j))||\nabla h(q_j)||<\dfrac{1}{j}\\
\nonumber 	\xi_1 (\sigma\circ h)(q_j)&>-\dfrac{1}{j}
 \end{align}
 
 then $ \mathcal{H}(h_\star)\ge H_2 $
 
 Thus $\mathcal{H}  $ is a constant function, $ \mathcal{H}'>0 $, so $ h $ is a constant function.

 (2)$\theta\ge 0,\rho'\le 0 $ is the same as (1) by replacing operator $ \zeta_1 $ by $ -\zeta_1 $.
\end{proof}

\vspace{2.12mm}
\subsection*{Acknowledgment}
I am grateful to Professor Hengyu Zhou(Chongqing University) for many useful discussions
and suggestions.

\end{CJK}
\end{document}